\def \be{\begin{equation}}
\def \eq{\end{equation}}
\def \bee{\begin{eqnarray}}
\def \eqq{\end{eqnarray}}

\def \zb{\bar{z}}

\def \R {{\cal R}}
\def \V {{\cal V}}
\def \W {{\cal W}}
\def \U {{\cal U}}
\def \A {{\cal A}}

\def \nn {\nonumber}

\documentstyle[12pt]{article}

\begin{document}
\begin{titlepage}
\begin{center}
July 12, 1995
        \hfill  LBL-37506\\
          \hfill    UCB-PTH-95/25\\
\vskip .3in

{\large \bf The Braided Quantum 2-Sphere}
\footnote{This work was supported in part by the Director, Office of
Energy Research, Office of High Energy and Nuclear Physics, Division of
High Energy Physics of the U.S. Department of Energy under Contract
DE-AC03-76SF00098 and in part by the National Science Foundation under
grant PHY-90-21139.}
\vskip .2in
{Chong-Sun Chu, Pei-Ming Ho and Bruno Zumino}
\footnote{email address: cschu@physics.berkeley.edu,
pmho@physics.berkeley.edu}

\vskip .2in
{\em
Department of Physics \\University of California \\ and\\
   Theoretical Physics Group\\
    Lawrence Berkeley Laboratory\\
      University of California\\
    Berkeley, CA 94720}
\end{center}

\vskip .2in
\begin{abstract}
In a recent paper the quantum 2-sphere $S^2_q$ was described as a quantum
complex manifold. Here we consider several copies of $S^2_q$ and derive their
braiding commutation relations. The braiding is extended to the differential
and to the integral calculus on the spheres. A quantum analogue of the
classical
anharmonic ratio of four points on the sphere is given, which is invariant
under
the coaction of $SU_q(2)$.

\end{abstract}

\end{titlepage}
\renewcommand{\thepage}{\roman{page}}
\setcounter{page}{2}
\mbox{ }

\vskip 1in

\begin{center}
{\bf Disclaimer}
\end{center}
\vskip .2in
\begin{scriptsize}
\begin{quotation}
This document was prepared as an account of work sponsored by the United
States Government. While this document is believed to contain correct
information, neither the United States Government nor any agency
thereof, nor The Regents of the University of California, nor any of their
employees, makes any warranty, express or implied, or assumes any legal
liability or responsibility for the accuracy, completeness, or usefulness
of any information, apparatus, product, or process disclosed, or represents
that its use would not infringe privately owned rights.  Reference herein
to any specific commercial products process, or service by its trade name,
trademark, manufacturer, or otherwise, does not necessarily constitute or
imply its endorsement, recommendation, or favoring by the United States
Government or any agency thereof, or The Regents of the University of
California.  The views and opinions of authors expressed herein do not
necessarily state or reflect those of the United States Government or any
agency thereof of The Regents of the University of California and shall
not be used for advertising or product endorsement purposes.
\end{quotation}
\end{scriptsize}

\vskip 2in

\begin{center}
\begin{small}
{\it Lawrence Berkeley Laboratory is an equal opportunity employer.}
\end{small}
\end{center}
\newpage
\renewcommand{\thepage}{\arabic{page}}
\setcounter{page}{1}

\section{Introduction}
In a recent paper\cite{CHZ} we considered the quantum sphere $S^2_q$ in the
sense of Podle\'{s} \cite{P1,P2,P3,P4}
for the particular value $c=0$ of his parameter $c$. We showed that
the quantum sphere can be described in terms of stereographic
variables $z$ and $\zb$ and we developed the differential and
integral calculus on $S^2_q$.
The entire formalism is covariant under the coaction
of the quantum group $SU_q(2)$ which is described by
quantum fractional transformations
on the variables $z$ and $\zb$, see Eq.(\ref{transf-z}) below.
The quantum sphere appears as a q-deformation of the
classical sphere considered as a complex K\"{a}hler manifold.
A $*$-conjugation antihomomorphism exists.

In this letter we consider the braiding of several copies of $S_q^2$.
There exists a
general formulation \cite{M1} for obtaining the braiding
of quantum spaces in terms of
the universal R matrix of the quantum group which
coacts on the quantum space. Using
this formulation the braiding commutation relations are obtained in Section 3
directly for the variables $z$ and $\zb$. (An alternative derivation of the
same
braiding relations proceeds by first computing the braiding of two copies of
the
complex quantum plane on which $SU_q(2)$ coacts and then using the expressions
of the stereographic variables $z$ and $\zb$ in terms of the coordinates
$x,y$ of the quantum plane
$$ z=x y^{-1}, \quad \zb=\bar{y}^{-1} \bar{x}. )$$

The braiding can be extended to the differentials $dz$ and $d\zb$. In Section 4
the
braiding property of the $SU_q(2)$ invariant integral on the sphere is given.
It is shown that it can be used to compute the integral.

The braiding of two quantum spheres is not symmetric with respect to
the exchange of
the two spheres.
It can be extended to the case of an arbitrary number of spheres
given in a certain order. For the case of four spheres one can construct a
quantum analogue of the classical anharmonic
ratio (cross ratio) of four points on a sphere.
This quantity, which belongs to the braided algebra of the four spheres, is
invariant under the coaction of $SU_q(2)$ as realized by the quantum fractional
transformations on the stereographic variables. It commutes with its
$*$-conjugate. This is described in Section 5.
The existence of the invariant anharmonic ratio seems very remarkable.

The results of this letter and those of \cite{CHZ} can
be extended to suitably defined
quantum deformations of
complex projective space and complex Grassmann manifolds.
These appear therefore as examples of quantum
complex K\"{a}hler manifolds. We are
planning to describe all this in a forthcoming publication.

\section{Braiding for Quantum Group Comodules}
Let $\A$ be the algebra of functions on a  quantum group and $\V$  an
algebra on which $\A$ coacts on the left:
\bee \Delta_L: \V &\rightarrow& \A \otimes \V \nn \\
               v & \mapsto & v^{(1')} \otimes v^{(2)}, \nn
\eqq
where we have used the Sweedler-like notation for $\Delta_L(v)$.

Let $\W$ be another left $\A$-comodule algebra,
\bee
   \Delta_L: \W &\rightarrow& \A \otimes \W \nn \\
               w & \mapsto & w^{(1')} \otimes w^{(2)}. \nn
\eqq
It is known\cite{M1} that one can
put $\V$ and $\W$ into a single left $\A$-comodule algebra
with the multiplication between
elements of $\V$ and $\W$ given by
\be \label{bcr} vw = \R(w^{(1')},v^{(1')}) w^{(2)} v^{(2)}. \eq
Here $\R \in \U \otimes \U$ is the universal R matrix for the quantum
enveloping
algebra $\U$ dual to $\A$ (with respect to the pairing $<\cdot, \cdot>$)
and
\[ \R(a,b) =<\R, a\otimes b>. \]
It satisfies:
\bee
  &\label{R1} \R(fg, h)=\R(f,h_{(1)}) \R(g,h_{(2)}), \\
  &\label{R2} \R(f, gh)=\R(f_{(1)},h) \R(f_{(2)},g), \\
  &\label{R3} \R(1,f) = \R(f,1) = \epsilon(f).
\eqq

For $\A=SU_q(2)$, it is
\[ \R(T^i_j, T^k_l) = q^{-1/2} \hat{R}^{ki}_{jl},
\]
where $\hat{R}$ is the $GL_q(2)$ R-matrix.
For example,
\bee &\R(a, T)= \pmatrix{q^{1/2}& 0\cr 0&q^{-1/2}},
  &\R(b, T)= \pmatrix{0& 0\cr 0&0}, \nn \\
     &\R(d, T)= \pmatrix{q^{-1/2}& 0\cr 0&q^{1/2}},
  &\R(c, T)= \pmatrix{0& \lambda q^{-1/2}\cr 0&0}, \nn
\eqq
where $\lambda=q-q^{-1}$.

\section{The Braided Sphere}
We first recall that the complex sphere $S_q^2$ is described by
coordinates $z, \zb$ which obey the commutation relation \cite{CHZ}
\be \label {zz}
 z \zb =q^{-2} \zb z +q^{-2}-1
\eq
and the $*$-structure $z^*=\zb$.
It is covariant under the fractional transformation,
with $\pmatrix{a&b\cr c&d} \in SU_q(2)$,
\be \label{transf-z}
z \rightarrow (a z +b)(c z+d)^{-1}, \quad
\zb \rightarrow -(c-d \zb)(a-b \zb)^{-1},
\eq
where $a,b,c$ and $d$ commute with $z$ and $\zb$.

One can extend $SU_q(2)$ by introducing $a^{-1}, d^{-1}$
satisfying
\bee
   &a a^{-1}= a^{-1} a =1, \quad d d^{-1}= d^{-1} d =1, \nn \\
   &\epsilon(a^{-1}) =\epsilon(d^{-1})=1, \nn \\
   &a^{-1 *} =d^{-1},\quad d^{-1 *}=a^{-1} \nn
\eqq
and
coproduct
\bee \Delta(d^{-1})&=& (c \otimes b + d \otimes d)^{-1} \nn \\
                    &=& \sum_{n=0}^{\infty} (-1)^n
            (d^{-1}c)^n d^{-1} \otimes (d^{-1}b)^n d^{-1}. \nn
\eqq
The transformation for $z$ and the braided copy $z'$ can then be written as
\[
z \rightarrow \sum_{n=0}^{n=\infty} f_n z^n,
\quad z' \rightarrow \sum_{n=0}^{n=\infty} f_n z'^n,
\]
where
\bee f_0&=&bd^{-1}, \nn \\
    f_n&=&(-1)^{n-1}d^{-2} (cd^{-1})^{n-1}, \quad n \geq 1 \nn \eqq
and Eq.(\ref{bcr}) gives
\bee \label{b-zz'} zz' &=& \sum_{n,m =0}^{\infty} \R(f_m,f_n) z'^m z^n.  \eqq

To calculate $\R(f_m, f_n)$,  we notice that, for example,
Eqs.(\ref{R1})-(\ref{R3}) give
\[ \R(a, a_{(1)}) \R(a^{-1}, a_{(2)}) =\R(a,a) \R(a^{-1}, a) =1 \]
and so
\bee \R(a^{-1}, T)= \pmatrix{q^{-1/2}& 0\cr 0&q^{1/2}},
     &\R(T,a^{-1})= \pmatrix{q^{-1/2}& 0\cr 0&q^{1/2}}, \nn
\eqq
\bee \R(d^{-1}, T)= \pmatrix{q^{1/2}& 0\cr 0&q^{-1/2}},
     &\R(T,d^{-1})= \pmatrix{q^{1/2}& 0\cr 0&q^{-1/2}}. \nn
\eqq

It is not hard to prove that for any functions $f,g$ of
$a^{\pm 1}, b,c, d^{\pm 1}$,
\bee &\R(bf, g) =0, &\R(f,cg)=0, \nn \eqq
\[ \R(cf,g)=0 ,\quad \mbox{if $g$ has no $b$}, \]
\[ \R(g,bf)=0, \quad \mbox{if $g$ has no $c$} \]
and
\be \label{Rd1} \R(d^{\pm 1},f(a,b,c,d))= f(q^{\mp 1/2},0,0,q^{\pm 1/2}) \eq
\be \label{Rd2} \R(f(a,b,c,d), d^{\pm 1} ) = f(q^{\mp 1/2},0,0,q^{\pm 1/2}) \eq
together with Eqs.(\ref{Rd1}),(\ref{Rd2})
with $d^{\pm 1}$ replaced by $a^{\mp 1}$.

One then gets
\bee \R(f_1, f_1)&=&q^2, \nn \\
     \R(f_2,f_0) &=&-\lambda q, \nn \\
     \R(f_m, f_n) &=& 0, \quad \mbox{all other $n,m$}. \nn
\eqq
Therefore
\be z z' =q^2 z' z -\lambda q z'^2. \label{zz'} \eq
Similarly,
\[ \zb \rightarrow \sum_{n=0}^{\infty} g_n \zb^n, \]
\bee g_0&=&-c a^{-1} , \nn \\
      g_n&=&q^{-2(n-1)} (ba^{-1})^{n-1} a^{-2}, \quad n \geq 1 \nn \eqq
and
\bee \R(g_1, f_1)&=&q^{-2}, \nn \\
     \R(g_0,f_0) &=&-\lambda q^{-1}, \nn \\
     \R(g_m, f_n) &=& 0, \quad \mbox{all other $n,m$}. \nn
\eqq
Therefore
\be \label{b-zzb'} z \zb' =
\sum_{n,m=0}^{\infty} \R(g_m, f_n) \zb^{'m} z^n   \eq
gives
\be z \zb' =q^{-2} \zb' z -\lambda q^{-1}. \label{zzb'} \eq

The differential calculus can also be defined on the braided spheres
by imposing the Leibniz rule on the exterior derivatives $d$ and $d'$
so that $d'$ acts on $z'$ and $\zb'$ in the same way
$d$ acts on $z$ and $\zb$, and
\bee
   &dz'=z'd, & d\zb'=\zb' d, \nn \\
   &d'z=zd', & d'\zb=\zb d' \nn
\eqq
together with
$$ d d' = -d' d. $$
Then Eqs.(\ref{zz'}),(\ref{zzb'}) and their $*$-involution
will imply  commutation relations between functions and forms
of different copies of the sphere.
As a consequence, the area element of the second copy
$K'=dz'd\zb'(1+\zb' z')^{-2}$ is central in the whole
braided algebra, while $K=dzd\zb(1+\zb z)^{-2}$ is only
central in the original copy $(z,\zb)$.

\section{ Remarks on the Integration }

For symmetry with respect to  the $*$-involution of the algebra
the braiding order of $z, \zb, z'$ and $\zb'$
has to be $z<(z',\zb')<\zb$
after we have determined $z<z'$ and $z<\zb'$
as in Eqs.(\ref{b-zz'}) and (\ref{b-zzb'}).
Because $z'$ and $\zb'$ are always on the same side
of the variables of their braided copy, $z$ and $\zb$,
the integration on $z',\zb'$,
has the following property:

if
\[ f(z',\zb')g(z,\zb)=\sum_{i}g_{i}(z,\zb)f_{i}(z',\zb'), \]
then
\be <f(z',\zb')>g(z,\zb)=\sum_{i}g_{i}(z,\zb)<f_{i}(z',\zb')>,
\label{fg} \eq
where $<\cdot>$ is the invariant integral on $S_{q}^{2}$.
This can be shown by using the invariance of the integral under the
$SU_q(2)$ coaction.
However,
\[ f(z',\zb')<g(z,\zb)>\neq\sum_{i}<g_{i}(z,\zb)>f_{i}(z',\zb'). \]

The above property (\ref{fg}) can be used to derive explicit integral rules.
For example, consider the case of $f(z',\zb')=\zb'\rho'^{-n}$,
where $\rho'=1+\zb' z'$ and $g(z,\zb)=z$.
Since
\[ \zb'\rho'^{-n}z = q^2 z\zb'\rho'^{-n}+
q^{1-2n}\lambda([n+1]_{q}-[n]_{q}\rho')\rho'^{-n},\quad n\geq 0, \]
where $[n]_{q}=\frac{q^{2n}-1}{q^{2}-1}$,
using Eq.(\ref{fg}) and $<\zb'\rho'^{-n}>=0$ we get
the recursion relation:
\be [n+1]_{q}<\rho'^{-n}>=[n]_{q}<\rho'^{-(n-1)}>, \quad n\geq 1.
\label{rr} \eq

This adds one more to two  other  different ways of computing the
invariant integrals.
One way is to impose the invariance condition directly \cite{CHZ}:
\[ <\chi f>=0, \]
for any generator $\chi$ of the $SU_q(2)$ quantum Lie algebra.
The other way is to use the cyclic property\footnote{
Similar cyclic properties have been found by H. Steinacker\cite{Ste1}
for integrals over higher dimensional quantum spheres in quantum
Euclidean space.}
\[ <f(z,\zb)g(z,\zb)>=<g(z,\zb)f(q^{-2}z,q^{2}\zb)>. \]
They all give the same recursion relation (\ref{rr}).

\section{ Anharmonic Ratios }

Let us first review the classical case.
The coordinates $x,y$ on a plane transform as
\be \label{transf-xy}
    \left(\begin{array}{c}
            x \\
            y
         \end{array}\right)\rightarrow
   \left(\begin{array}{cc}
            a & b \\
            c & d
         \end{array}\right)
   \left(\begin{array}{c}
            x \\
            y
         \end{array}\right)
\eq
by an $SU(2)$ matrix
$T=\pmatrix{a&b\cr c&d}$.
The determinant-like object
$xy'-yx'$ defined for $x,y$ together with
the coordinates on a second plane $x',y'$
is invariant under the $SU(2)$ transformation.
On each plane we define $z=x/y$ so that
\[ z-z'=y^{-1}(xy'-yx')y'^{-1}. \]
It now follows that with $x_i,y_i$ for $i=1,2,3,4$
as coordinates on four copies of the planes,
\bee & &(z_2-z_1)(z_2-z_4)^{-1}(z_3-z_4)(z_3-z_1)^{-1} \nn \\
&=&(x_1 y_2-y_1 x_2)(x_4 y_2-y_4 x_2)^{-1}
(x_4 y_3-y_4 x_3)(x_1 y_3-y_1 x_3)^{-1} \nn
\eqq
is invariant because all the factors $y_i^{-1}$ cancel
and only the invariant parts $(x_i y_j-y_i x_j)$ survive.
Therefore the anharmonic ratio is invariant.

Permuting the indices in the above expression
we may get other anharmonic ratios,
but they are all functions of the one above.
For example,
\[
(z_2-z_3)(z_2-z_4)^{-1}(z_1-z_4)(z_3-z_1)^{-1}=
(z_2-z_1)(z_2-z_4)^{-1}(z_3-z_4)(z_3-z_1)^{-1}-1.
\]

The $SU_q(2)$ covariant quantum plane obeys
$$ xy=qyx $$
which is covariant under  the  transformation (\ref{transf-xy})
with $T$ now being an $SU_q(2)$ matrix.
Braided quantum planes can be introduced by using Eq.(\ref{bcr}).
Let  $\V$ be the $i$-th copy and $\W$ be the $j$-th one, then
we have for $i<j$,
\bee &x_i y_j =qy_j x_i +q \lambda x_j y_i, \nn \\
     &x_i x_j =q^2 x_j x_i, \nn \\
     &y_i y_j =q^2 y_j y_i, \nn \\
     &y_i x_j =q x_j y_i. \nn
\eqq
In the deformed case
we have to be more careful about the ordering.
Let the deformed determinant-like object be
\[ (ij):=x_i y_j-q y_i x_j, \]
which is invariant under the $SU_q(2)$ transformation,
and let
\[ [ij] := z_i - z_j = q^{-1} y_i^{-1} (ij) y_j^{-1}, \]
where $z_i := x_i y_i^{-1}$.

Using the relations
\[ y_i(ij)=q (ij)y_i, \]
\[ (ij)y_j=q y_j(ij) \]
for $i<j$ and
\[ y_i(jk)=q^3 (jk)y_i, \]
\[ (ij)y_k=q^3 y_k(ij) \]
for $i<j<k$,
we can see that, for example,
\[ A:=[12][24]^{-1}[34][13]^{-1} \]
is again invariant.
Similarly,
$B:=[12][23]^{-1}[34][14]^{-1}$
as well as a number of others are invariant.

To find out whether  these invariants are
independent of one another,
we now discuss the algebra of the $[ij]$'s.

Because $[ij]=[ik]+[kj]$ and $[ij]=-[ji]$
the algebra of $[ij]$ for $i,j=1,2,3,4$ is generated
by only three elements $[12],[23],[34]$.
It is easy to prove that
\[
   [ij][kl]=q^{2}[kl][ij],
\]
if $i<j\leq k<l$.

It follows that we have
\be
   [ij][ik][jk]=q^4[jk][ik][ij] \label{YB}
\eq
for $i<j<k$, and
\[
   [12][34]+[14][23]=[12][24]+[24][23].
\]
Using  these relations we can check the dependency between the
different anharmonic
ratios.
For example, let $C:=[13][23]^{-1}[24][14]^{-1}$,
and $D:=[14][13]^{-1}[23][24]^{-1}$, both invariant,
then
\bee
   B^{-1}AC &=& [14][34]^{-1}[23][24]^{-1}([34][23]^{-1}[24])[14]^{-1} \nn \\
            &=& [14][34]^{-1}[23][24]^{-1}([24][23]^{-1}[34])[14]^{-1} \nn \\
            &=& 1, \nn
\eqq
where we used the relation $[34][23]^{-1}[24]=[24][23]^{-1}[34]$
which follows Eq.(\ref{YB}),
and
\bee
   q^{2}B-D^{-1} &=& ([12][34][23]^{-1}-[24][23]^{-1}[13])[14]^{-1} \nn \\
                  &=& ([12][34][23]^{-1}-[24][23]^{-1}([12]+[23]))
                      [14]^{-1} \nn \\
                  &=& ([12][34]-[12][24]-[24][23])[23]^{-1}[14]^{-1} \nn \\
                  &=& (-[14][23])[23]^{-1}[14]^{-1} \nn \\
                  &=& -1. \nn
\eqq

In this manner it can be checked that
all products of four terms
$[ij]$, $[kl]$, $[mn]^{-1}$, $[pr]^{-1}$
in arbitrary order, which are invariant, are  functions of
only one invariant, say, $A$.
Namely, all invariants are related and
just like in the classical case,
there is only one independent anharmonic ratio.
It can be checked that the anharmonic ratio commutes
with all the $\zb_i$'s and so commutes with its $*$-complex
conjugate, which is also an invariant.

\section{Acknowledgements}
We are very grateful to Chryssomalis Chryssomalakos, Michael Schlieker
and Harold Steinacker for numerous very helpful discussions.
This work was supported in part by the Director, Office of
Energy Research, Office of High Energy and Nuclear Physics, Division of
High Energy Physics of the U.S. Department of Energy under Contract
DE-AC03-76SF00098 and in part by the National Science Foundation under
grant PHY-90-21139.

\end{document}